\documentstyle[twoside]{article}
\include{graphicx}
\title{\bf Asymptotic expansion of the modified exponential integral involving the Mittag-Leffler function}
\author{\sc R. B. Paris\footnote{E-mail address:\ \ {\tt r.paris@abertay.ac.uk}}\\
\\
{\em Division of Computing and Mathematics,}\\
{\em Abertay University, Dundee DD1 1HG, UK}\\
}

\begin{document}
\newcommand{\bee}{\begin{equation}}
\newcommand{\ee}{\end{equation}}
\def\f#1#2{\mbox{${\textstyle \frac{#1}{#2}}$}}
\def\dfrac#1#2{\displaystyle{\frac{#1}{#2}}}
\newcommand{\fr}{\frac{1}{2}}
\newcommand{\fs}{\f{1}{2}}
\newcommand{\g}{\Gamma}
\newcommand{\al}{\alpha}
\newcommand{\bb}{\beta}
\newcommand{\om}{\omega}
\newcommand{\br}{\biggr}
\newcommand{\bl}{\biggl}
\newcommand{\ra}{\rightarrow}
\renewcommand{\topfraction}{0.9}
\renewcommand{\bottomfraction}{0.9}
\renewcommand{\textfraction}{0.05}
\newcommand{\mcol}{\multicolumn}
\newcommand{\gtwid}{\raisebox{-.8ex}{\mbox{$\stackrel{\textstyle >}{\sim}$}}}
\newcommand{\ltwid}{\raisebox{-.8ex}{\mbox{$\stackrel{\textstyle <}{\sim}$}}}
\date{}
\maketitle
\pagestyle{myheadings}
\markboth{\hfill {\it R.B. Paris} \hfill}
{\hfill {\it Asymptotics of a modified exponential integral} \hfill}
\begin{abstract} 
We consider the asymptotic expansion of the generalised exponential integral involving the Mittag-Leffler function introduced recently by Mainardi and Masina [{\it Fract. Calc. Appl. Anal.} {\bf 21} (2018) 1156--1169]. We extend the definition of this function using the two-parameter Mittag-Leffler function. The expansions of the similarly extended  sine and cosine integrals are also discussed.
Numerical examples are presented to illustrate the accuracy of each type of expansion obtained. 
\vspace{0.4cm}

\noindent {\bf MSC:} 30E15, 30E20, 33E20, 34E05 
\vspace{0.3cm}

\noindent {\bf Keywords:} asymptotic expansions, exponential integral, Mittag-Leffler function, sine and cosine inegrals\\
\end{abstract}

\vspace{0.2cm}

\noindent $\,$\hrulefill $\,$

\vspace{0.2cm}

\begin{center}
{\bf 1. \  Introduction}
\end{center}
\setcounter{section}{1}
\setcounter{equation}{0}
\renewcommand{\theequation}{\arabic{section}.\arabic{equation}}
The complementary exponential integral $\mbox{Ein}(z)$ is defined by
\bee\label{e10}
\mbox{Ein}(z)=\int_0^z \frac{1-e^{-t}}{t}\,dt=\sum_{n=1}^\infty \frac{(-)^{n-1}z^n}{n n!}\ \ \ (z\in\bf{C})
\ee
and is an entire function.
Its connection with the classical exponential integral ${\cal E}_1(z)=\int_z^\infty t^{-1}e^{-t}\,dt$, valid in the cut plane $|\arg\,z|<\pi$, is \cite[p.~150]{DLMF}
\bee\label{e11}
\mbox{Ein}(z)=\log\,z+\gamma+{\cal E}_1(z),
\ee
where $\gamma=0.5772156\dots$ is the Euler-Mascheroni constant.

In a recent paper, Mainardi and Masina \cite{MM} proposed an extension of $\mbox{Ein}(z)$ by replacing the exponential function in (\ref{e10}) by the one-parameter Mittag-Leffler function
\[E_\alpha(z)=\sum_{n=0}^\infty\frac{z^n}{\g(\alpha n+1)}\qquad (z\in{\bf C},\ \ \alpha>0),\]
which generalises the exponential function $e^z$. They introduced the function for any $\alpha>0$ in the cut plane $|\arg\,z|<\pi$
\bee\label{e12a}
\mbox{Ein}_\alpha(z)=\int_0^z\frac{1-E_\alpha(-t^\alpha)}{t^\alpha}\,dt=\sum_{n=0}^\infty \frac{(-)^n z^{\alpha n+1}}{(\alpha n+1) \g(\alpha n+\al+1)},
\ee
which when $\alpha=1$ reduces to the function $\mbox{Ein}(z)$. A physical application of this function for $0\leq\alpha\leq 1$ arises in the study of the creep features of a linear viscoelastic model; see \cite{MMS} for details. An analogous extension of the generalised sine and cosine integrals was also considered in \cite{MM}.
Plots of all these functions for $\alpha\in[0,1]$ were given.

Here we consider a slightly more general version of (\ref{e12a}) based on the two-parameter Mittag-Leffler function given by
\[E_{\alpha,\beta}(z)=\sum_{n=0}^\infty\frac{z^n}{\g(\alpha n+\beta)}\qquad(z\in{\bf C},\ \alpha>0),\]
where $\bb$ will be taken to be real.
Then the extended complementary exponential integral we shall consider is
\[\mbox{Ein}_{\alpha,\beta}(z)=\int_0^z\frac{1-E_{\alpha,\beta}(-t^\alpha)}{t^\alpha}\,dt=\sum_{n=1}^\infty \frac{(-)^{n-1}}{\g(\alpha n+\beta)}\int_0^z t^{\alpha n-\alpha}dt\]
\bee\label{e21}
=z\sum_{n=0}^\infty\frac{(-)^{n} z^{\alpha n}}{(\alpha n+1) \g(\alpha n+\alpha+\beta)}\hspace{2cm}
\ee
upon replacement of $n-1$ by $n$ in the last summation. When $\beta=1$ this reduces to (\ref{e12a}) so that $\mbox{Ein}_{\alpha,1}(z)=\mbox{Ein}_\alpha(z)$. 

The asymptotic expansion of this function will be obtained for large complex $z$ with the parameters $\al$, $\bb$ held fixed.
We achieve this by consideration of the asymptotics of a related function using the theory developed for integral functions of hypergeometric type as discussed, for example, in \cite[\S2.3]{PK}. An interesting feature of the expansion of $\mbox{Ein}_{\al,\bb}(x)$ for $x\to+\infty$ when $\al\in(0,1]$ is the appearance of a logarithmic term
whenever $\al=1,\fs, \f{1}{3}, \ldots\ $.
Similar expansions are obtained for the extended sine and cosine integrals in Section 4. The paper concludes with the presentation of some numerical results that demonstrate the accuracy of the different expansions obtained.

\vspace{0.6cm}

\begin{center}
{\bf 2. \ The asymptotic expansion of a related function for $|z|\to\infty$}
\end{center}
\setcounter{section}{2}
\setcounter{equation}{0}
\renewcommand{\theequation}{\arabic{section}.\arabic{equation}}
To determine the asymptotic expansion of $\mbox{Ein}_{\alpha,\beta}(z)$ for large complex $z$ with the parameters $\alpha$ and $\beta$ held fixed, we shall find it convenient to consider the related function defined by
\bee\label{e22}
F(\chi):=\sum_{n=0}^\infty \frac{\chi^n}{(\alpha n+\gamma) \g(\alpha n+\al+\bb)}=\sum_{n=0}^\infty g(n) \frac{\chi^n}{n!}
\qquad(\chi\in{\bf C}),
\ee
where 
\[g(n)=\frac{\g(n+1)}{(\al n+\gamma)\g(\al n+\al+\bb)}=\frac{\g(\al n+\gamma) \g(n+1)}{\g(\al n+\gamma+1)\g(\al n+\al+\bb)}~.\]
The parameter $\gamma>0$, but will be chosen to have two specific values in Sections 3 and 4; namely, $\gamma=1$ and $\gamma=1+\al$. 
It will be shown that the asymptotic expansion of $F(\chi)$ consists of an algebraic and an exponential expansion valid in different sectors of the complex $\chi$-plane.

The function $F(\chi)$ in (\ref{e22}) is a case of the Wright function 
\begin{equation}\label{e23}
{}_p\Psi_q(\chi)=\sum_{n=0}^\infty \frac{\prod_{r=1}^p\Gamma(\alpha_rn+a_r)}{\prod_{r=1}^q\Gamma(\beta_rn+b_r)}\,\frac{\chi^n}{n!}, 
\end{equation}
corresponding to $p=q=2$. In (\ref{e23}) the parameters $\alpha_r$  and 
$\beta_r$ are real and positive and $a_r$ and $b_r$ are
arbitrary complex numbers. We also assume that the $\alpha_r$ and $a_r$ are subject to 
the restriction
\[\alpha_rn+a_r\neq 0, -1, -2, \ldots \qquad (n=0, 1, 2, \ldots\ ;\, 1\leq r \leq p)\]
so that no gamma function in the numerator in (\ref{e23}) is singular.
The parameters associated\footnote{Empty sums and products are to be interpreted as zero and unity, respectively.} with ${}_p\Psi_q(\chi)$ are given by
\[\kappa=1+\sum_{r=1}^q\beta_r-\sum_{r=1}^p\alpha_r, \qquad 
h=\prod_{r=1}^p\alpha_r^{\alpha_r}\prod_{r=1}^q\beta_r^{-\beta_r},\]
\begin{equation}\label{e12}
\vartheta=\sum_{r=1}^pa_r-\sum_{r=1}^qb_r+\f{1}{2}(q-p),\qquad \vartheta'=1-\vartheta.
\end{equation}

The algebraic expansion of $F(\chi)$ is obtained from the Mellin-Barnes integral representation \cite[p.~56]{PK}
\[F(\chi)=\frac{1}{2\pi i}\int_{c-\infty i}^{c+\infty i} \frac{\g(-s) \g(1+s)(\chi e^{\mp\pi i})^s}{(\al s+\gamma) \g(\al s+\al+\bb)}\,ds,\qquad |\arg (-\chi)|<\pi(1-\fs\al),\]
where, with $-\gamma/\al<c<0$, the integration path lies to the left of the poles of $\g(-s)$ at $s=0, 1, 2, \ldots$ but to the right of the poles at $s=-\gamma/\al$ and $s=-k-1$, $k=0, 1, 2, \ldots\,$. The upper or lower sign is taken according as $\arg\,\chi>0$ or $\arg\,\chi<0$, respectively. It is seen that when $\al=\gamma/m$, $m=1, 2, \ldots$ the pole at $s=-m$ is double and its residue must be evaluated accordingly. Displacement of the integration path to the left when $0<\al<2$ and evaluation of the residues  then produces the algebraic expansion $H(\chi e^{\mp\pi i})$, where
\bee\label{e25}
H(\chi)=\left\{\begin{array}{ll}\!\!\!\dfrac{\pi/\al}{\sin \gamma\pi/\al}\,\dfrac{\chi^{-\gamma/\al}}{\g(\al+\bb-\gamma)}+\sum_{k=0}^\infty \dfrac{(-)^k \chi^{-k-1}}{(\gamma-\al(k+1)) \g(\bb-\al k)} & (\al\neq \dfrac{\gamma}{m})\\
\\
\!\!\!\dfrac{(-)^{m-1} \chi^{-m}}{\g(\al+\bb-\gamma)}\bl\{\frac{m}{\gamma}\log\,\chi-\psi(\al+\bb-\gamma)\br\}+\mathop{\sum_{k=0}^\infty}_{\scriptstyle k\neq m-1} \dfrac{(-)^k \chi^{-k-1}}{(\gamma-\al(k+1)) \g(\bb-\al k)} & (\al=\dfrac{\gamma}{m}),\end{array}\right.
\ee
and $\psi$ denotes the logarithmic derivative of the gamma function.

The exponential expansion associated with ${}_p\Psi_q(\chi)$ is given by \cite[p.~299]{PFCH}, \cite[p.~57]{PK}
\begin{equation}\label{e22c}
{\cal E}(\chi):=X^\vartheta e^X\sum_{j=0}^\infty A_jX^{-j}, \qquad X=\kappa (h\chi)^{1/\kappa},
\end{equation}
where the coefficients $A_j$ are those appearing in the inverse factorial expansion  
\begin{equation}\label{e22a}
\frac{1}{\g(1+s)}\,\frac{\prod_{r=1}^p\Gamma(\alpha_rn+a_r)}{\prod_{r=1}^q\Gamma(\beta_rn+b_r)}=\kappa A_0 (h\kappa^\kappa)^s\bl\{\sum_{j=0}^{M-1}\frac{c_j}{\Gamma(\kappa s+\vartheta'+j)}
+\frac{\rho_M(s)}{\Gamma(\kappa s+\vartheta'+M)}\br\}
\end{equation}
with $c_0=1$.
Here $M$ is a positive integer and $\rho_M(s)=O(1)$ for $|s|\ra\infty$ in $|\arg\,s|<\pi$.
The constant $A_0$ is specified by
\[
A_0=(2\pi)^{\frac{1}{2}(p-q)}\kappa^{-\frac{1}{2}-\vartheta}\prod_{r=1}^p
\alpha_r^{a_r-\frac{1}{2}}\prod_{r=1}^q\beta_r^{\frac{1}{2}-b_r}.
\]
The coefficients $c_j$ are independent of $s$ and depend only on the parameters $p$, $q$, $\alpha_r$, 
$\beta_r$, $a_r$ and $b_r$.

For the function $F(\chi)$, we have 
\[\kappa=\al,\quad h=\al^{-\al},\quad \vartheta=-\al-\bb, \quad A_0=\al^{-1}.\]
We are in the fortunate position that the normalised coefficients $c_j$ in this case can be determined explicitly  as $c_j=(\al+\bb-\gamma)_j$. This follows from the well-known (convergent) expansion given in \cite{F}, \cite[p.~41]{PK}
\bee\label{e24c}
\frac{1}{(\al s+\gamma)\g(\al s+\al+\bb)}=\sum_{j=0}^\infty \frac{(\al+\bb-\gamma)_j}{\g(\al s+\vartheta'+j)}\qquad (\Re (s)>-\gamma/\al),
\ee
to which, in the case of $F(\chi)$, the ratio of gamma functions appearing on the left-hand side of (\ref{e22a})  reduces.
Then, with $X=\chi^{1/\al}$ we have from (\ref{e22c}) the exponential expansion associated with $F(\chi)$ given by
\bee\label{e24}
{\cal E}(\chi)=\frac{1}{\al} \chi^{\vartheta/\al} \exp\,[\chi^{1/\al}]\,\sum_{j=0}^\infty\, (\al+\bb-\gamma)_j \,\chi^{-j/\al}.
\ee

From \cite[pp.~57--58]{PK}, we then obtain the asymptotic expansion  for $|\chi|\to\infty$ when $0<\al<2$
\bee\label{e26}
F(\chi)\sim\left\{\begin{array}{ll} {\cal E}(\chi)+H(\chi e^{\mp\pi i}) & |\arg\,\chi|<\fs\pi\al \\
\\
H(\chi e^{\mp\pi i}) & |\arg (-\chi)|<\pi(1-\fs\al)\end{array}\right.\ee
and, when $\al=2$,
\bee\label{e27}
F(\chi) \sim {\cal E}(\chi)+ {\cal E}(\chi e^{\mp2\pi i})+ H(\chi e^{\mp\pi i})\qquad |\arg\,\chi|\leq\pi.
\ee
The upper and lower signs are chosen according as $\arg\,\chi>0$ or $\arg\,\chi<0$, respectively. It may be noted that the expansions ${\cal E}(\chi e^{\mp2\pi i})$ in (\ref{e27}) only become significant in the neighbourhood of $\arg\,\chi=\pm\pi$. 
When $\al>2$, the expansion of $F(\chi)$ is exponentially large for all values of $\arg\,\chi$ (see \cite[p.~58]{PK}) and accordingly we omit this case as it is unlikely to be of physical interest.
\vspace{0.2cm}

\noindent{\bf Remark 2.1}\ \ The exponential expansion ${\cal E}(\chi)$ in (\ref{e26}) continues to hold beyond the sector $|\arg\,\chi|<\fs\pi\al$, where it becomes exponentially small in the sectors $\pi\al\leq|\arg\,\chi|<\fs\pi\al$ when $0<\al\leq 1$. The rays $\arg\,\chi=\pm\pi\al$ are Stokes lines, where ${\cal E}(\chi)$ is maximally subdominant relative to the algebraic expansion $H(\chi e^{\mp\pi i})$. On these rays, ${\cal E}(\chi)$ undergoes a Stokes phenomenon, where the exponentially small expansion ``switches off'' in a smooth manner as $|\arg\,\chi|$ increases \cite[\S 2.11(iv)]{DLMF}, with its value to leading order given by $\fs{\cal E}(\chi)$; see \cite{P13} for a more detailed discussion of this point in the context of the confluent hypergeometric functions. We do not consider exponentially small contributions to $F(\chi)$ here, except to briefly mention in Section 3 the situation pertaining to the case $\al=1$.

\vspace{0.6cm}

\begin{center}
{\bf 3. \ The asymptotic expansion of $\mbox{Ein}_{\al,\bb}(z)$ for $|z|\to\infty$}
\end{center}
\setcounter{section}{3}
\setcounter{equation}{0}
\renewcommand{\theequation}{\arabic{section}.\arabic{equation}}
The asymptotic expansion of $\mbox{Ein}_{\al,\bb}(z)$ defined in (\ref{e21}) can now be constructed from that of $F(\chi)$ with the parameter $\gamma=1$. It is sufficient, for real $\al$, $\bb$, to consider $0\leq\arg\,z\leq\pi$, since the expansion when $\arg\,z<0$ is given by the conjugate value. With $\chi=-z^\al=e^{-\pi i} z^\al$, the exponentially large sector $|\arg\,\chi|<\fs\pi\al$ becomes $|-\pi+\al \arg\,z|<\fs\pi\al$; that is
\bee\label{e27a}
\theta_0<\arg\,z < \theta_0+\pi,\qquad\theta_0:=\frac{\pi}{2\al}(2-\al).
\ee
On the boundaries of this sector the exponential expansion is of an oscillatory character. When $0<\al<\f{2}{3}$, we note that the exponentially large sector (\ref{e27a}) lies outside the sector of interest $0\leq\arg\,z\leq\pi$.

We define the algebraic and exponential asymptotic expansions
\bee\label{e28}
H_{\alpha,\beta}(z)=\left\{\begin{array}{ll}\!\!\! \dfrac{\pi/\al}{\sin (\pi/\al)\,\g(\al+\bb-1)}+\sum_{k=0}^\infty\dfrac{(-)^k z^{1-\al(k+1)}}{(1-\al(k+1)) \g(\bb-\al k)} & (\al\neq m^{-1})\\
\\
\!\!\!\dfrac{(-)^{m-1}}{\g(\al+\bb-1)}\{\log\,z-\psi(\al+\bb-1)\}+\mathop{\sum_{k=0}^\infty}_{\scriptstyle k\neq m-1}\dfrac{(-)^k z^{1-\al (k+1)}}{(1-\al(k+1)) \g(\bb-\al k)} & (\al=m^{-1}),\end{array}\right.
\ee
where $m=1, 2, \ldots\,$, and
\bee\label{e29}
{\cal E}_{\al,\bb}(z)=\frac{(e^{-\pi i/\al}z)^{\vartheta}}{\al}\,\exp\,[e^{-\pi i/\al} z] \sum_{j=0}^\infty (\al+\bb-1)_j (e^{-\pi i/\al}z)^{-j}
\ee
where we recall that $\vartheta=-\al-\bb$.
Then the following result holds:
\vspace{0.3cm}

\noindent{\bf Theorem 1.}\ \ {\it Let $m$ be a positive integer, with $\al>0$ and $\bb$ real and $\theta_0=\pi(2-\al)/(2\al)$. Then the following expansions hold for $|z|\to\infty$}
\bee\label{e210a}
\mbox{Ein}_{\al,\bb}(z)\sim H_{\al,\bb}(z)\qquad (0\leq\arg\,z\leq\pi)
\ee
{\it when $0<\al<\f{2}{3}$, and}
\bee\label{e10b}
\mbox{Ein}_{\al,\bb}(z)\sim\left\{\begin{array}{ll}H_{\al,\bb}(z) & (0\leq\arg\,z<\theta_0)\\
\\
z{\cal E}_{\al,\bb}(z)+H_{\al,\bb}(z) & (\theta_0\leq\arg\,z\leq\pi)\end{array}\right.
\ee
{\it when $\f{2}{3}\leq\al<2$. Finally, when $\al=2$ we have} $\mbox{Ein}_{2,\bb}(-z)=-\mbox{Ein}_{2,\bb}(z)$ {\it and it is therefore sufficient to consider $0\leq\arg\,z\leq\fs\pi$. Then, from (\ref{e27}), we obtain the expansion when $\al=2$}
\bee\label{e210c}
\mbox{Ein}_{2,\bb}(z)\sim z\{{\cal E}_{2,\bb}(z)+ {\cal E}_{2,\bb}(ze^{\pi i})\}+H_{2,\bb}(z)\qquad(0\leq\arg\,z\leq\fs\pi).
\ee
\vspace{0.1cm}

\noindent We note from Theorem 1 that when $z\to-\infty$ the value of $\mbox{Ein}_{\al,\bb}(z)$ is, in general, complex-valued.

In the case of main physical interest, when $z=x>0$ is a real variable, we have the following expansion:
\vspace{0.3cm}

\noindent{\bf Theorem 2.}\ \ {\it When $z=x$ $(>0)$ we have from Theorem 1 the expansions}
\bee\label{e211a}
\mbox{Ein}_{\al,\bb}(x)\sim H_{\al,\bb}(x) 
\ee
{\it for $0<\al<2$, and from (\ref{e29}) and (\ref{e210c}) when $\al=2$}
\bee\label{e211b}
\mbox{Ein}_{2,\bb}(x)\sim H_{2,\bb}(x)-x^{-1-\bb}\sum_{j=0}^\infty\frac{(1+\bb)_j}{x^j}\,\cos\,[x-\fs\pi(\bb+j)]
\ee
{\it as $x\to+\infty$.}
\vspace{0.2cm}

\noindent It is worth noting that a logarithmic term is present in the asymptotic expansion of $\mbox{Ein}_{\al,\beta}(x)$ whenever $\al=1, \fs, \f{1}{3}, \ldots\ $.
\vspace{0.2cm}

\noindent{\bf 3.1}\ \ {\it The case $\al=1$}
\vspace{0.1cm}

The special case $\al=1$ deserves further consideration. From (\ref{e28}) and (\ref{e211a}) we obtain the expansion
\bee\label{e212a}
\mbox{Ein}_{1,\bb}(x)\sim \frac{1}{\g(\bb)}\{\log\,x-\psi(\bb)\}-\sum_{k=1}^\infty\frac{(-x)^{-k}}{k\, \g(\bb-k)}\qquad (x\to+\infty).
\ee
If $\bb=1$, the asymptotic sum in (\ref{e212a}) vanishes and
\bee\label{e212}
\mbox{Ein}_{1,1}(x)\sim \log\,x+\gamma
\ee
for large $x$. But we have the exact evaluation (compare (\ref{e11}))
\[\mbox{Ein}_{1,1}(x)=x\sum_{n=0}^\infty\frac{(-x)^n}{(n+1)^2 n!}=\log\,x+\gamma+{\cal E}_1(x)\]
\bee\label{e213}
\hspace{2.4cm} \sim \log\,x+\gamma+\frac{e^{-x}}{x}\sum_{j=0}^\infty \frac{(-)^j j!}{x^j}\qquad (x\to+\infty)
\ee
by \cite[(6.12.1)]{DLMF}.
The additional asymptotic sum appearing in (\ref{e213}) is exponentially small as $x\to+\infty$ and is consequently not accounted for in the result (\ref{e212}).

From Remark 2.1, it is seen that there are Stokes lines at $\arg\,z=\pm\pi(1-\al)$, which coalesce on the positive real axis when $\al=1$. In the sense of increasing $\arg\,z$ in the neighbourhood of the positive real axis, the exponential expansion ${\cal E}_{1,\bb}(z)$ is in the process of {\it switching on} across $\arg\,z=\pi(1-\al)$ and 
$\overline{{\cal E}_{1,\bb}(z)}$ (where the bar denotes the complex conjugate) is in the process of {\it switching off\/} across $\arg\,z=-\pi(1-\al)$. When $\al=1$, this produces the exponential contribution
\[\fs x\{{\cal E}_{1,\bb}(x)+\overline{{\cal E}_{1,\bb}(x)}\}=\frac{e^{-x}}{x^\bb}\,\cos \pi\bb\,\sum_{j=0}^\infty \frac{(-)^{j+1} (\bb)_j}{x^j}\]
for large $x$. Thus, the more accurate version of (\ref{e212a}) should read
\bee\label{e214}
\mbox{Ein}_{1,\bb}(x)\sim\frac{1}{\g(\bb)}\{\log\,x-\psi(\bb)\}-\sum_{k=1}^\infty\frac{(-x)^{-k}}{k\, \g(\bb-k)}
-\frac{e^{-x}}{x^\bb}\,\cos \pi\bb\,\sum_{j=0}^\infty \frac{(-)^{j} (\bb)_j}{x^j}
\ee
as $x\to+\infty$. When $\bb=1$, this correctly reduces to (\ref{e213}).

When $\bb=2$, we have \cite{M}
\[\mbox{Ein}_{1,2}(x)=x\sum_{n=0}^\infty \frac{(-x)^n}{(n+1)^2(n+2) n!}=\log\,x-\psi(2)+\frac{1}{x}+{\cal E}_1(x)-\frac{e^{-x}}{x}\]
\[\hspace{0.4cm}\sim \log\,x-\psi(2)+\frac{1}{x}+\frac{e^{-x}}{x}\sum_{j=1}^\infty \frac{(-)^j j!}{x^j}\qquad (x\to+\infty).\]
This can be seen also to agree with (\ref{e214}) after a little rearrangement.

\vspace{0.6cm}

\begin{center}
{\bf 4. \ The generalised sine and cosine integrals}
\end{center}
\setcounter{section}{4}
\setcounter{equation}{0}
\renewcommand{\theequation}{\arabic{section}.\arabic{equation}}
The sine and cosine integrals are defined by \cite[\S 6.2]{DLMF}
\[\mbox{Si}(z)=\int_0^z\frac{\sin t}{t}\,dt,\qquad \mbox{Cin}(z)=\int_0^z\frac{1-\cos t}{t}\,dt.\]
Mainardi and Masina \cite{MM} generalised these definitions by replacing the trigonometric functions by
\[\sin_\al(t)=t^\al E_{2\al,\al+\bb}(-t^{2\al})=\sum_{n=0}^\infty \frac{(-)^n t^{(2n+1)\al}}{\g(2n\al+\al+\bb)},\quad
\cos_\al(t)=E_{2\al,\bb}(-t^{2\al})=\sum_{n=0}^\infty \frac{(-)^n t^{2n\al}}{\g(2n\al+\bb)}\]
with $\bb=1$ to produce
\bee\label{e41}
\left\{\begin{array}{l}
\mbox{Sin}_\al(z)={\displaystyle\int_0^z\dfrac{\sin_\al(t)}{t^\al}\,dt=\sum_{n=0}^\infty \dfrac{(-)^n z^{2n\al+1}}{(2n\al+1) \g(2n\al+\al+1)}}\\
\\
\mbox{Cin}_\al(z)={\displaystyle\int_0^\infty\dfrac{1-\cos_\al(t)}{t^\al}\,dt=\sum_{n=0}^\infty \dfrac{(-)^n z^{2n\al+\al+1}}{(2n\al+\al+1) \g(2n\al+2\al+1)}}~.\end{array}\right.
\ee
Here we extend the definitions (\ref{e41}) by including the additional parameter $\bb\in{\bf R}$ in the Mittag-Leffler functions and consider the functions
\bee\label{e42}
\left\{\begin{array}{l}
\mbox{Sin}_{\al,\bb}(z)={\displaystyle z\sum_{n=0}^\infty \dfrac{(-)^n z^{2n\al}}{(2n\al+1) \g(2n\al+\al+\bb)}}\\
\\
\mbox{Cin}_{\al,\bb}(z)={\displaystyle z^{1+\al}\sum_{n=0}^\infty \dfrac{(-)^n z^{2n\al}}{(2n\al+\al+1) \g(2n\al+2\al+\bb)}}~.\end{array}\right.
\ee

The asymptotics of $\mbox{Sin}_{\al,\bb}(z)$ and $\mbox{Cin}_{\al,\bb}(z)$ can be deduced from the results in Section 2. However, here we restrict ourselves to determining the asymptotic expansion of these functions for large $|z|$ in a sector enclosing the positive real $z$-axis, where for $0<\al<1$ they only have an algebraic-type expansion.
We observe in passing that
\bee\label{e40}
\mbox{Sin}_{\al,\bb}(z)=\mbox{Ein}_{2\al,\bb-\al}(z).
\ee

Comparison of the series expansion for $\mbox{Sin}_{\al,\bb}(z)$ with $F(\chi)$ in Section 2, with the substitutions $\al\to 2\al$, $\bb\to\bb-\al$ and $\gamma=1$ (or from the above identity combined with Theorems 1 and 2), produces the following expansion:
\vspace{0.3cm}

\noindent{\bf Theorem 3.}\ \ {\it For $m=1, 2, \ldots$ and\/ $0<\al<1$ we have the algebraic expansions}
\[\mbox{Sin}_{\al,\bb}(z)\hspace{13cm}\]
\bee\label{e43}
\sim\left\{\begin{array}{ll}
\dfrac{\pi/(2\al)}{\sin (\pi/(2\al)) \g(\al\!+\!\bb\!-\!1)}+\sum_{k=0}^\infty \dfrac{(-)^k z^{1-2\al(k+1)}}{(1-2\al(k+1)) \g(\bb-(2k+1)\al)} & (\al\neq(2m)^{-1})\\
\\
\dfrac{(-)^{m-1}}{\g(\al\!+\!\bb\!-\!1)}\{\log\,z\!-\!\psi(\al\!+\!\bb\!-\!1)\}+ \mathop{\sum_{k=0}^\infty}_{\scriptstyle k\neq m-1} \dfrac{(-)^k z^{1-2\al(k+1)}}{(1-2\al(k+1)) \g(\bb-(2k+1)\al)} & (\al=(2m)^{-1})\end{array}\right.
\ee
{\it as $|z|\to\infty$ in the sector $|\arg\,z|<\pi(1-\al)/(2\al)$.}
\vspace{0.3cm}

A similar treatment for $\mbox{Cin}_{\al,\bb}(z)$ shows that with the substitutions $\al\to2\al$, $\bb\to\bb$ and $\gamma=1+\al$ we obtain the following expansion:
\vspace{0.3cm}

\noindent{\bf Theorem 4.}\ \ {\it For $m=1, 2, \ldots$ and\/ $0<\al<1$ we have the algebraic expansions}
\[\mbox{Cin}_{\al,\bb}(z)\hspace{13cm}\]
\bee\label{e44}
\sim\left\{\begin{array}{ll}
\dfrac{\pi/(2\al)}{\cos (\pi/(2\al)) \g(\al\!+\!\bb\!-\!1)}+\sum_{k=0}^\infty \dfrac{(-)^k z^{1-(2k+1)\al }}{(1-(2k+1)\al) \g(\bb-2\al k)} & (\al\neq(2m-1)^{-1})\\
\\
\dfrac{(-)^{m-1}}{\g(\al\!+\!\bb\!-\!1)}\{\log\,z\!-\!\psi(\al\!+\!\bb\!-\!1)\}+ \mathop{\sum_{k=0}^\infty}_{\scriptstyle k\neq m-1} \dfrac{(-)^k z^{1-(2k+1)\al}}{(1-(2k+1)\al) \g(\bb-2\al k)} & (\al=(2m-1)^{-1})\end{array}\right.
\ee
{\it as $|z|\to\infty$ in the sector $|\arg\,z|<\pi(1-\al)/(2\al)$.}
\vspace{0.3cm} 

\noindent The expansions of $\mbox{Sin}_{\al,\bb}(x)$ and $\mbox{Cin}_{\al,\bb}(x)$ as $x\to+\infty$ when $0<\al<1$ follow immediately from Theorems 3 and 4.

As $x\to+\infty$ when $\al=1$, the exponentially oscillatory contribution to $\mbox{Sin}_{1,\bb}(x)$ can be obtained directly from (\ref{e211b}) together with (\ref{e40}). In the case of $\mbox{Cin}_{1,\bb}(x)$, we obtain from (\ref{e22c}) with $\kappa=2$, $h=\f{1}{4}$, $\vartheta=-2-\bb$, $X=\chi^{1/2}$ and $A_0=\fs$ the exponential expansion
\[{\cal E}(\chi)=\frac{1}{2}\chi^{\vartheta/2} \exp\,[\chi^{1/2}] \sum_{j=0}^\infty c_j \chi^{-j/2}, \qquad \chi=e^{-\pi i}x^2,\]
with the coefficients $c_j=(\bb)_j$. Then the exponential contribution to $\mbox{Cin}_{1,\bb}(x)$ is 
\[x^2\{{\cal E}(\chi)+{\cal E}(\chi e^{\pi i})\}=-x^{-\bb}\sum_{j=0}^\infty \frac{(\bb)_j}{x^j}\,\cos\,[x-\fs\pi(\bb+j)]\qquad (x\to+\infty).\]

 Collecting together these results we finally obtain the following theorem.
\vspace{0.3cm}

\noindent{\bf Theorem 5.}\ \ {\it When $\al=1$ and $\bb$ is real the following expansions hold:}
\[\mbox{Sin}_{1,\bb}(x)\sim\frac{\pi}{2\g(\bb)}-\sum_{k=0}^\infty\frac{(-)^k x^{-2k-1}}{(2k+1) \g(\bb-1-2k)}\hspace{4cm}\]
\bee\label{e45}
\hspace{6cm}+x^{-\bb}\sum_{j=0}^\infty \frac{(\bb)_j}{x^j}\,\sin \,[x-\fs\pi(\bb+j)]
\ee
and
\[\mbox{Cin}_{1,\bb}(x) \sim \frac{1}{\g(\bb)}\,\{\log\,x-\psi(\bb)\}-\sum_{k=1}^\infty \frac{(-)^k x^{-2k}}{2k \g(\bb-2k)}\hspace{3.3cm}\]
\bee\label{e46}
\hspace{6cm}-x^{-\bb} \sum_{j=0}^\infty\frac{(\bb)_j}{x^j}\,\cos\,[x-\fs\pi(\bb+j)]
\ee
{\it as $x\to+\infty$.}
\vspace{0.3cm}

\noindent When $\bb>0$, it is seen that $\mbox{Sin}_{1,\bb}(x)$ approaches the constant value $\pi/(2\g(\bb))$ whereas $\mbox{Cin}_{1,\bb}(x)$ grows logarithmically like $\log\,(x) /\g(\bb)$ as $x\to+\infty$. 

\vspace{0.6cm}

\begin{center}
{\bf 5. \ Numerical results}
\end{center}
\setcounter{section}{5}
\setcounter{equation}{0}
\renewcommand{\theequation}{\arabic{section}.\arabic{equation}}
In this section we present numerical results confirming the accuracy of the various expansions obtained in this paper. In all cases we have employed optimal truncation (that is truncation at, or near, the least term in modulus) of the algebraic and (when appropriate) the exponential expansions. 

We first present results in the physically interesting case of $0<\al\leq 1$ and $\bb=1$ considered in \cite{MM}.
Table 1 shows the values\footnote{In the tables we write the values as $x(y)$ instead of $x\times 10^y$.} of the absolute relative error in the computation of $\mbox{Ein}_{\al,1}(x)$ from the asymptotic expansions in Theorem 2
for several values of $x$ and different $\al$ in the extended range $0<\al\leq 2$. The expansion for $0<\al<2$ is given by the algebraic expansion in (\ref{e211a}); this contains a logarithmic term for the values $\al=\f{1}{4}, \fs, 1$. The progressive loss of accuracy when $\al>1$ can be attributed to the presence of the approaching exponentially large sector, whose lower boundary is, from (\ref{e27a}), given by $\theta_0=\pi(2-\al)/(2\al)$.
In the final case $\al=2$, the accuracy is seen to suddenly increase considerably. This is due to the inclusion of the (oscillatory) exponential contribution, which from (\ref{e211b}), takes the form
\[\mbox{Ein}_{2,1}(x)\sim \frac{1}{2}\pi-\frac{1}{x}-\sum_{j=1}^\infty\frac{j!}{x^{j+1}}\,\cos (x-\fs\pi j)\qquad(x\to+\infty).\] 
In Fig. 1 we show some plots of $\mbox{Ein}_{\al,1}(x)$ for values of $\al$ in the range $0<\al\leq 1$. In Fig. 2 the asymptotic approximations for two values of $\al$ are shown compared with the corresponding curves of $\mbox{Ein}_{\al,1}(x)$.

\begin{table}[h]
\caption{\footnotesize{The absolute relative error in the computation of $\mbox{Ein}_{\al,\bb}(x)$ from Theorem 2
for different values of $\al$ and $x$ when $\bb=1$.}}
\begin{center}
\begin{tabular}{|l|ccccc|}
\hline
&&&&&\\[-0.25cm]
\mcol{1}{|c|}{$x$} & \mcol{1}{c}{$\al=0.25$} &\mcol{1}{c}{$\al=0.40$} & \mcol{1}{c}{$\al=0.50$} & \mcol{1}{c}{$\al=0.75$} & \mcol{1}{c|}{$\al=1.00$}\\
[.1cm]\hline
&&&&&\\[-0.25cm]
5  & $1.602(-04)$ & $1.678(-05)$ & $2.012(-04)$ & $2.115(-04)$ & $5.249(-04)$ \\
10 & $5.733(-07)$ & $1.735(-07)$ & $4.413(-07)$ & $2.339(-07)$ & $1.442(-06)$ \\
20 & $3.680(-11)$ & $3.031(-11)$ & $6.526(-12)$ & $9.362(-12)$ & $2.753(-11)$ \\
30 & $1.808(-16)$ & $9.384(-16)$ & $1.543(-16)$ & $1.337(-16)$ & $7.595(-16)$ \\
[.1cm]\hline
&&&&&\\[-0.25cm]
\mcol{1}{|c|}{$x$} & \mcol{1}{c}{$\al=1.20$} &\mcol{1}{c}{$\al=1.40$} & \mcol{1}{c}{$\al=1.60$} & \mcol{1}{c}{$\al=1.80$} & \mcol{1}{c|}{$\al=2.00$}\\
[.1cm]\hline
&&&&&\\[-0.25cm]
5  &  $1.121(-03)$ & $1.301(-04)$ & $5.279(-03)$ & $1.407(-02)$ & $1.550(-03)$ \\
10 &  $4.345(-06)$ & $3.168(-05)$ & $2.103(-04)$ & $1.536(-04)$ & $2.849(-06)$ \\
20 &  $2.147(-10)$ & $2.277(-08)$ & $1.671(-06)$ & $4.751(-05)$ & $4.926(-10)$ \\
30 &  $4.388(-14)$ & $2.363(-11)$ & $2.125(-08)$ & $6.216(-06)$ & $1.613(-14)$ \\
[.1cm]\hline
\end{tabular}
\end{center}
\end{table}
\begin{figure}[th]
	\begin{center}	{\tiny($a$)}\includegraphics[width=0.5\textwidth]{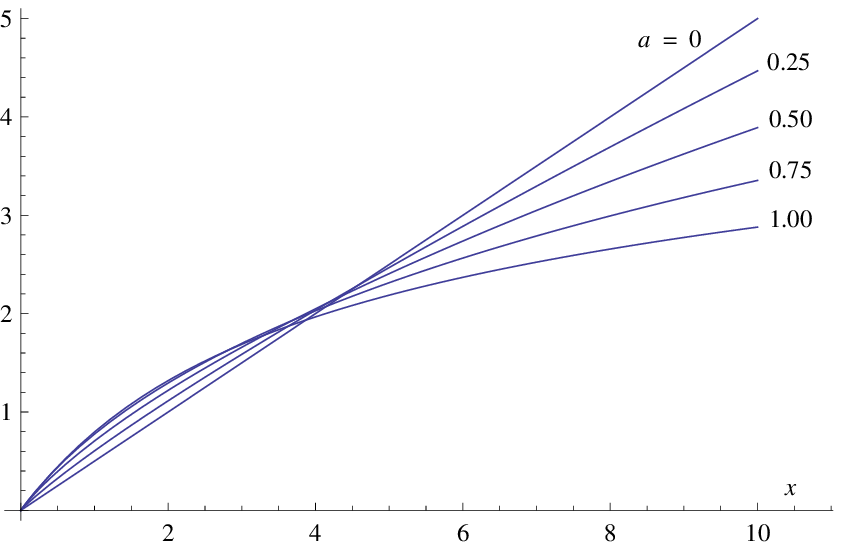}
\caption{\small{Plots of $\mbox{Ein}_{\al,1}(x)$ for different values of $\al$.
}}
	\end{center}
\end{figure}
\begin{figure}[th]
	\begin{center}	{\tiny($a$)}\includegraphics[width=0.375\textwidth]{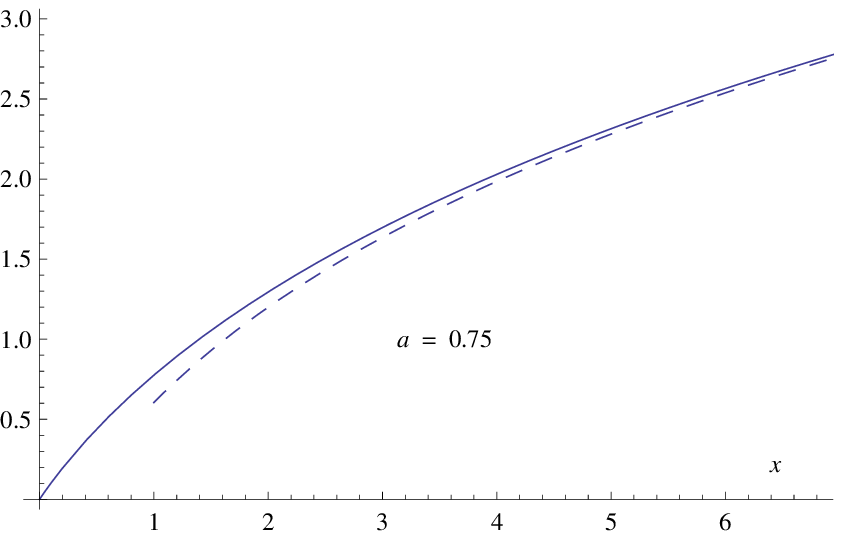}
	\qquad
	{\tiny($b$)}\includegraphics[width=0.375\textwidth]{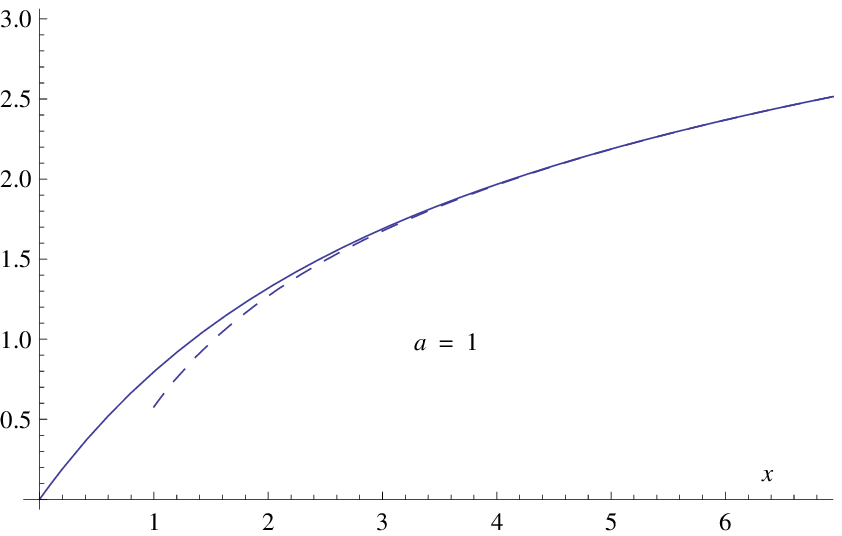}
\caption{\small{Plots of $\mbox{Ein}_{\al,1}(x)$ (solid curves) and the leading asymptotic approximation (dashed curves) for ($a$) $\al=0.75$ and ($b$) $\al=1$
}}
	\end{center}
\end{figure}

Table 2 shows the values of the absolute relative error in the computation of $\mbox{Ein}_{\al,\bb}(z)$ from the asymptotic expansions in Theorem 1 for complex $z$ for values of $\al$ in the range $0<\al\leq2$. It will noticed that there is a sudden reduction in the error when $\al=1$ and $\theta=\pi/4$.
In this case, the value of $\theta_0=\fs\pi$ and a more accurate treatment would include the exponentially small contribution $z{\cal E}_{\al,\bb}(z)$. When this term is included we find the absolute relative error equal to $6.935\times 10^{-11}$.

\begin{table}[h]
\caption{\footnotesize{The absolute relative error in the computation of $\mbox{Ein}_{\al,\bb}(z)$ from Theorem 1 
for different $\al$ and $\theta$ when $z=20e^{i\theta}$ and $\bb=1/3$.}}
\begin{center}
\begin{tabular}{|l|ccccc|}
\hline
&&&&&\\[-0.25cm]
\mcol{1}{|c|}{$\theta$} & \mcol{1}{c}{$\al=0.40$} &\mcol{1}{c}{$\al=0.50$} & \mcol{1}{c}{$\al=1.00$} & \mcol{1}{c}{$\al=1.50$} & \mcol{1}{c|}{$\al=2.00$}\\
[.1cm]\hline
&&&&&\\[-0.25cm]
0 &        $2.400(-08)$ & $5.494(-10)$ & $2.702(-10)$ & $1.572(-06)$ & $5.119(-10)$ \\
$\pi/4$ &  $2.553(-08)$ & $1.820(-09)$ & $1.142(-07)$ & $1.202(-08)$ & $8.204(-08)$ \\
$\pi/2$ &  $3.026(-08)$ & $4.057(-09)$ & $1.756(-10)$ & $2.021(-08)$ & $3.684(-07)$ \\
$3\pi/4$ & $3.897(-08)$ & $8.028(-09)$ & $1.423(-09)$ & $2.320(-07)$ & $8.204(-08)$ \\
$\pi$ &    $5.398(-08)$ & $1.617(-08)$ & $6.457(-09)$ & $3.005(-03)$ & $5.119(-10)$ \\
[.1cm]\hline
\end{tabular}
\end{center}
\end{table}

Finally, in Table 3 we present the error associated with the expansions of the generalised sine and cosine integrals $\mbox{Sin}_{\al,\bb}(x)$ and 
$\mbox{Cin}_{\al,\bb}(x)$ as $x\to+\infty$ given in Theorems 3--5. For $\mbox{Sin}_{\al,\bb}(x)$, the logarithmic expansion in (\ref{e43}) arises for $\al=\f{1}{4}$ and $\al=\fs$; for $\mbox{Cin}_{\al,\bb}(x)$ the logarithmic expansion in (\ref{e44}) arises for $\al=\f{1}{3}$. In Fig. 3 are shown plots\footnote{We remark that the plot of $\mbox{Cin}_{\al1}(x)$ in Fig.~3(b) differs from that shown in Fig.~4 of \cite{MM}.} of $\mbox{Sin}_{\al,1}(x)$ and $\mbox{Cin}_{\al,1}(x)$ for different $\al$ and in Fig. 4 the leading asymptotic approximations from the expansions in Theorem 5
%\[\mbox{Sin}_{1,1}(x)\sim \frac{1}{2}\pi-\sum_{j=0}^\infty \frac{j!}{x^{j+1}}\,\cos (x-\fs\pi j),\]
%\[\mbox{Cin}_{1,1}(x)\sim \log\,x+\gamma-\sum_{j=0}^\infty \frac{j!}{x^{j+1}}\,\sin (x-\fs\pi j)\]
are compared with the corresponding plots of these functions.
\begin{table}[th]
\caption{\footnotesize{The absolute relative error in the computation of $\mbox{Sin}_{\al,\bb}(x)$ and $\mbox{Cin}_{\al,\bb}(x)$ from Theorems 3--6 for different $\al$ and $x$ when $\bb=4/3$ .}}
\begin{center}
\begin{tabular}{|l|ccccc|}
\hline
&&&&&\\[-0.25cm]
\mcol{1}{|c|}{} & \mcol{5}{c|}{$\mbox{Sin}_{\al,\bb}(x)$}\\
\mcol{1}{|c|}{$x$} & \mcol{1}{c}{$\al=1/4$} &\mcol{1}{c}{$\al=1/3$} & \mcol{1}{c}{$\al=1/2$} & \mcol{1}{c}{$\al=2/3$} & \mcol{1}{c|}{$\al=1$}\\
[.1cm]\hline
&&&&&\\[-0.25cm]
10 &  $4.396(-07)$ & $1.394(-08)$ & $1.785(-06)$ & $3.410(-06)$ & $1.012(-05)$ \\
20 &  $3.213(-11)$ & $1.171(-13)$ & $3.920(-11)$ & $2.076(-08)$ & $3.094(-11)$ \\
25 &  $2.373(-13)$ & $3.792(-14)$ & $2.098(-13)$ & $4.437(-10)$ & $3.270(-12)$ \\
30 &  $1.879(-15)$ & $5.065(-15)$ & $1.172(-15)$ & $8.197(-12)$ & $8.010(-15)$ \\
[.1cm]\hline
\hline
&&&&&\\[-0.25cm]
\mcol{1}{|c|}{} & \mcol{5}{c|}{$\mbox{Cin}_{\al,\bb}(x)$}\\
\mcol{1}{|c|}{$x$} & \mcol{1}{c}{$\al=1/4$} &\mcol{1}{c}{$\al=1/3$} & \mcol{1}{c}{$\al=1/2$} & \mcol{1}{c}{$\al=2/3$} & \mcol{1}{c|}{$\al=1$}\\
[.1cm]\hline
&&&&&\\[-0.25cm]
10 &  $9.237(-08)$ & $3.787(-07)$ & $6.608(-07)$ & $2.270(-05)$ & $7.756(-06)$ \\
20 &  $1.293(-12)$ & $4.473(-12)$ & $1.090(-11)$ & $2.462(-10)$ & $2.576(-10)$ \\
25 &  $8.066(-14)$ & $2.334(-16)$ & $5.326(-14)$ & $6.881(-11)$ & $1.437(-12)$ \\
30 &  $1.160(-16)$ & $9.285(-17)$ & $2.764(-16)$ & $2.934(-12)$ & $7.716(-15)$ \\
[.1cm]\hline
\end{tabular}
\end{center}
\end{table}
\begin{figure}[th]
	\begin{center}	{\tiny($a$)}\includegraphics[width=0.375\textwidth]{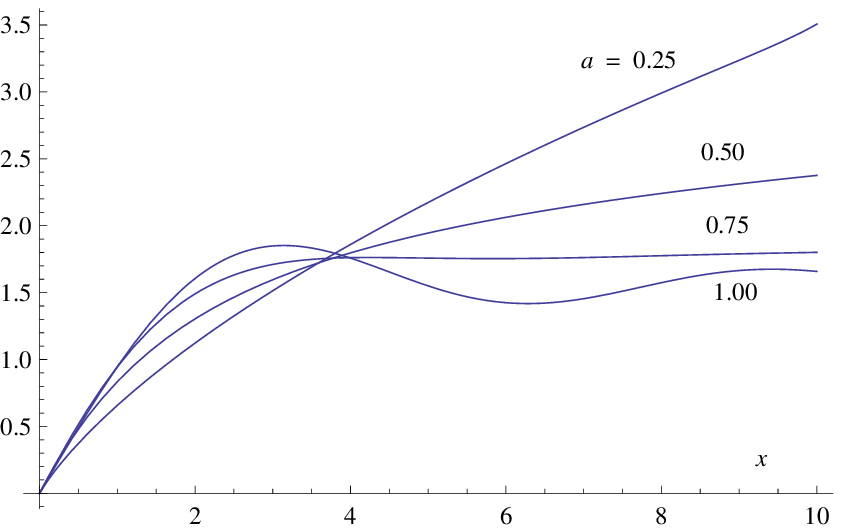}
	\qquad
	{\tiny($b$)}\includegraphics[width=0.375\textwidth]{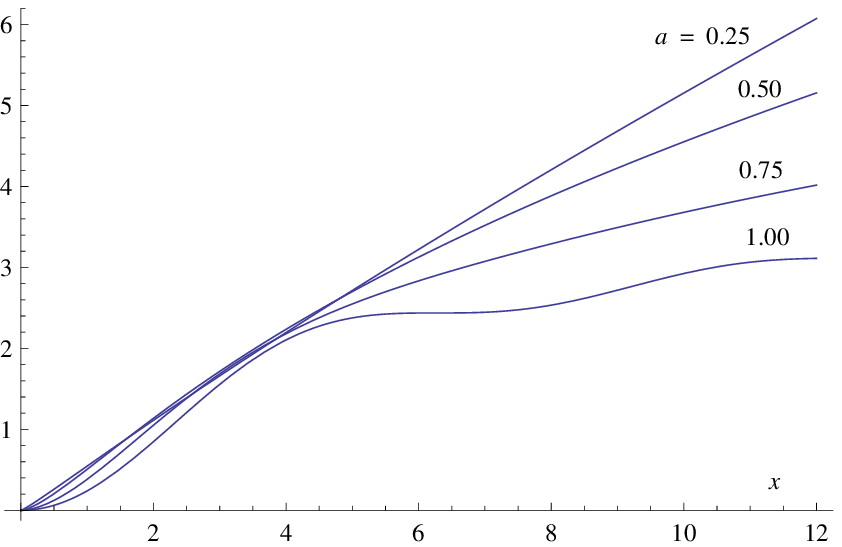}
	\vspace{0.4cm}
	
%	{\tiny($c$)}\includegraphics[width=0.375\textwidth]{TouchardFig1c}\qquad
%	{\tiny($c$)}\includegraphics[width=0.375\textwidth]{TouchardFig1d}
\caption{\small{Plots of the generalised sine and cosine integrals ($a$) $\mbox{Sin}_{\al,1}(x)$  and ($b$) $\mbox{Cin}_{\al,1}(x)$ for $\al=0.25, 0.50, 0.75,1$.
}}
	\end{center}
\end{figure}
\begin{figure}[th]
	\begin{center}	{\tiny($a$)}\includegraphics[width=0.375\textwidth]{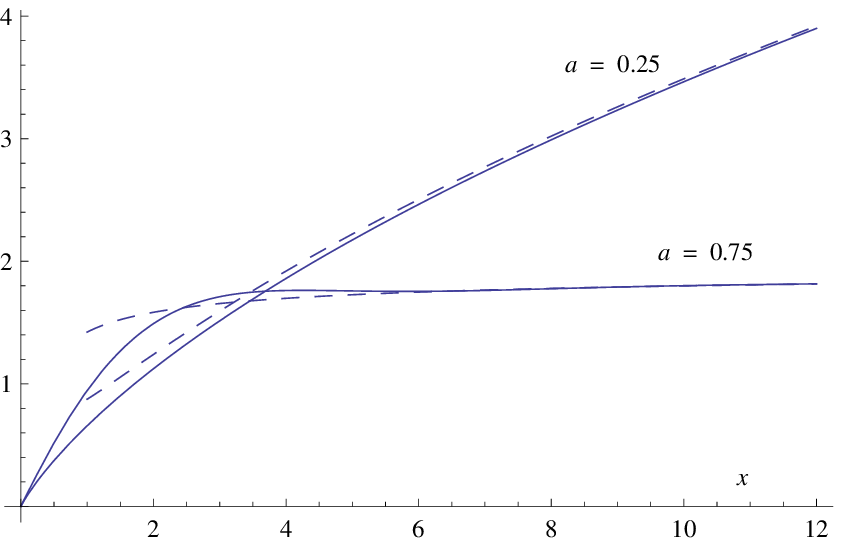}
	\qquad
	{\tiny($b$)}\includegraphics[width=0.375\textwidth]{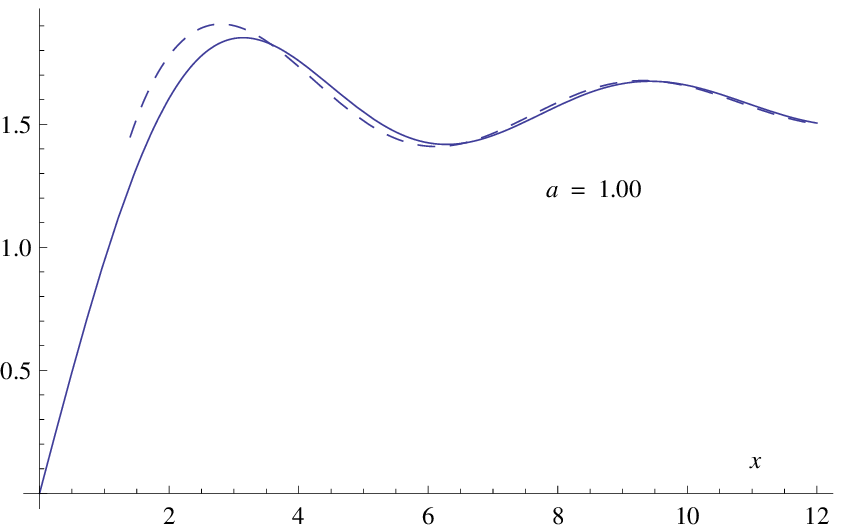}
	\vspace{0.4cm}
	
	{\tiny($c$)}\includegraphics[width=0.375\textwidth]{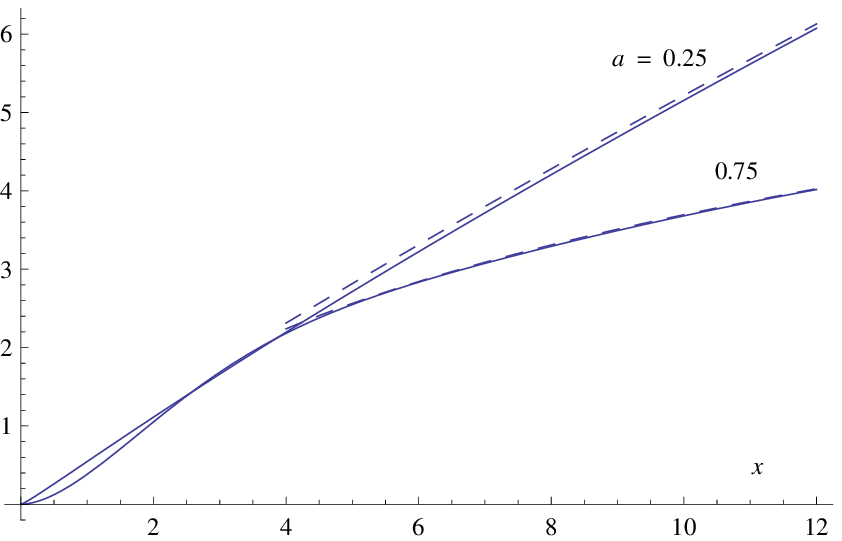}\qquad
	{\tiny($d$)}\includegraphics[width=0.375\textwidth]{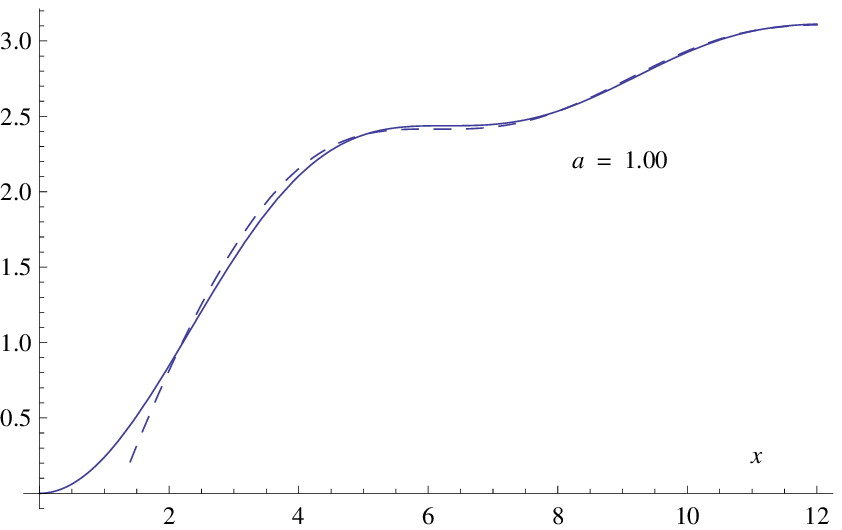}
\caption{\small{Plots of the generalised sine and cosine integrals (solid curves) and their leading asymptotic approximations (dashed curves) from Theorems 3, 4 and 5: ($a$) $\mbox{Sin}_{\al,1}(x)$ when $\al=0.25, 0.75$, ($b$) $\mbox{Sin}_{\al,1}(x)$ when $\al=1$, ($c$)  $\mbox{Cin}_{\al,1}(x)$ when $\al=0.25, 0.75$ and ($d$) $\mbox{Cin}_{\al,1}(x)$ when $\al=1$.
}}
	\end{center}
\end{figure}

\end{document}